\documentclass[12pt]{amsart}

\usepackage{amsthm}
\usepackage{amssymb}
\usepackage{amsmath}
\numberwithin{equation}{section}
\newcommand{\bea}{\begin{eqnarray}}
\newcommand{\eea}{\end{eqnarray}}
\newcommand{\be}{\begin{eqnarray*}}
\newcommand{\ee}{\end{eqnarray*}}

\begin{document}

\title[Gaussian binomials and sublattices]{Gaussian binomials and the number of sublattices}
\author[Yi Ming Zou]{Yi Ming Zou}
\address{Department of Mathematical Sciences, University of Wisconsin, Milwaukee, WI 53201, USA} \email{ymzou@uwm.edu}
\maketitle

\begin{abstract}
The purpose of this short communication is to make some observations on the connections between various existing formulas of counting the number of sublattices of a fixed index in an $n$-dimensional lattice and their connection with the Gaussian binomials.
\end{abstract}
\section{Existing formulas}
There are various ways of determining the number of sublattices of a fixed index in a lattice, they can be found in Cassels (1971), Baake (1997), and Gruber (1997). To determine the number $f_{n}(m)$ (notation as in Baake (1997)) of sublattices of index $m$ in an $n$-dimensional lattice is the same as to determine the number of subgroups of index $m$ in a free abelian group of rank $n$. A detailed discussion of this problem was included in Baake (1997), where a formula to compute $f_{n}(m)$:
\bea
f_{n}(m)=\sum_{d_{1}d_{2}\cdots d_{n}=m}d^{0}_{1}d^{1}_{2}\cdots d^{n-1}_{n},
\eea
and a formula to express the generating function $F_{n}(s)$ of $f_{n}(m)$ as a Dirichlet series ($\zeta(s)=\sum_{m=1}^{\infty}m^{-s}$ is the Riemann zeta function):
\bea
F_{n}(s) = \zeta(s)\zeta(s-1)\cdots\zeta(s-n+1),
\eea
are provided. These formulas imply the following recursion relation:
\be
f_{n}(m)=\sum_{d\mid m}df_{n-1}(d).
\ee
One can also use the results on pp.11-13 from Cassels (1971) to show that $f_{n}(m)$ is equal to the number of $n\times n$ matrices $(r_{ij})$ with integer entries satisfying the conditions (lower triangular matrices):
\bea
r_{ij}=0,\quad 1\le i<j\le n;\\\nonumber
r_{ii}>r_{ij}\ge 0, \quad 1\le j<i\le n;\\
r_{11}r_{22}\cdots r_{nn} = m.\nonumber
\eea
Gruber (1997) proved that if $m=p^{r_{1}}_{1}\cdots p^{r_{k}}_{k}$ is the prime factorization of $m$, then $f_{n}(m)$ can also be computed by the following formula:  
\bea
f_{n}(m) = \prod_{i=1}^{k}\prod_{j=1}^{r_{i}}\frac{p^{n+j-1}_{i}-1}{p^{j}_{i}-1}= \prod_{i=1}^{k}\prod_{j=1}^{n-1}\frac{p^{r_{i}+j}_{i}-1}{p^{j}_{i}-1}.
\eea
Although all these methods were mentioned in Gruber (1997), the connections among these existing methods have not been adequately explained. In the next section, we will make some observations on the connections among these methods as well as the connection with Gaussian binomials.
\section{Observations}
First, we observe that formula (1.1) can be derived from Cassels' result by noting that for an integer matrix satisfying condition (1.3), there are $r_{ii}$ choices for each of the elements below $r_{ii}$ at the $i$-th column, and therefore the number of these matrices for each decomposition $m=r_{11}r_{22}\cdots r_{nn}$ is $r_{11}^{n-1}r_{22}^{n-2}\cdots r_{nn}^{0}$. Summing over all decompositions, one gets (1.1).
\par
Then we observe that formula (1.4) can be derived from formula (1.2) by using the Gaussian binomials. Recall (Jantzen (1995) Ch. 0 or the online Wikipedia) that for integers $m, k\ge 0$, the Gaussian binomials ($q$-binomial coefficients) are defined by
\bea
\left[\begin{array}{c} m\\ k \end{array}\right]_{q} = \frac{[m]_{q}!}{[m-k]_{q}![k]_{q}!},
\eea   
where
\bea
[m]_{q}=\frac{1-q^{m}}{1-q},\quad [m]_{q}!=[1]_{q}[2]_{q}\cdots [m]_{q}.
\eea
These binomials satisfy 
\bea
\left[\begin{array}{c} m\\ k \end{array}\right]_{q}=\left[\begin{array}{c} m\\ m-k \end{array}\right]_{q}.
\eea
By using Gaussian binomials, one has the following formula
\bea
\prod_{k=0}^{n-1}\frac{1}{1-q^{k}t}=\sum_{k=0}^{\infty}\left[\begin{array}{c} n+k-1\\ k \end{array}\right]_{q}t^{k}.
\eea 
Now since
\be
F_{n}(s) &=& \prod_{i=0}^{n-1}\zeta(s-i)=\prod_{p}\prod_{k=0}^{n-1}\frac{1}{1-p^{-s+k}}\\
&=& \prod_{p}\sum_{k=0}^{\infty}\left[\begin{array}{c} n+k-1\\ k \end{array}\right]_{p}p^{-sk}\\
& =& \sum_{m=1}^{\infty}\frac{\prod_{i=1}^{r_{m}}\left[\begin{array}{c} n+k_{i}-1\\ k_{i} \end{array}\right]_{p_{i}}}{m^{s}},
\ee
where $m=p_{1}^{k_{1}}\cdots p_{r_{m}}^{k_{r_{m}}}$, we obtain the first formula in (1.4), and applying (2.3), we obtain the second formula in (1.4).
\section{Concluding remarks} 
There exist two approaches to the problem of counting the number of sublattices of a fixed index in the literature. The most detailed discussion is provided by Baake (1997). Alternatively, one can derive (1.1) from Cassels (1971), and then prove (1.2) by using arguments similar to that of Baake (1997). The connection between the number of sublattices of fixed indices and the Gaussian binomials is provided by the two product formulas in Gruber (1997).      
\par
\medskip

\end{document}